\documentclass[11pt]{article}
\usepackage[matrix,arrow,curve,cmtip]{xy}
\usepackage{amssymb}
\usepackage{latexsym}
\usepackage{theorem}
\usepackage[a4paper,width=160mm,height=230mm]{geometry}
\usepackage{titlesec}
\usepackage{graphicx}

\titleformat{\section}[hang]{\bf\Large}{\thesection.}{1ex}{}
\titleformat{\subsection}[hang]{\bfseries\normalsize}{\thesubsection.}{1ex}{}

\def\to{\mbox{$\xymatrix@1@C=5mm{\ar@{->}[r]&}$}}
\def\halfcirc{\begin{picture}(0,0)\put(0,2){\oval(4,4)[l]}\end{picture}}
\def\incl{\mbox{$\xymatrix@1@C=5mm{\ar@{->}[r]|<{\halfcirc}&}$}}
\def\tto{\mbox{$\xymatrix@1@C=5mm{\ar@{=>}[r]&}$}}
\def\distsign{\begin{picture}(0,0)\put(0,0){\circle{4}}\end{picture}}
\def\dist{\mbox{$\xymatrix@1@C=5mm{\ar@{->}[r]|{\distsign}&}$}}
\def\spansign{\begin{picture}(0,0)\put(0,-3){\line(0,1){6}}\end{picture}}
\def\span{\mbox{$\xymatrix@1@C=5mm{\ar@{->}[r]|{\spansign}&}$}}
\def\criblesign{\begin{picture}(0,0)\put(-1,-3){\line(0,1){6}}\put(1,-3){\line(0,1){6}}\end{picture}}
\def\crible{\mbox{$\xymatrix@1@C=5mm{\ar@{->}[r]|{\criblesign}&}$}}

\newtheorem{theorem}{Theorem}[section]
\newtheorem{lemma}[theorem]{Lemma}
\newtheorem{definition}[theorem]{Definition} 
\newtheorem{proposition}[theorem]{Proposition}
\newtheorem{corollary}[theorem]{Corollary}
{\theorembodyfont{\upshape}\newtheorem{example}[theorem]{Example}}
\newcommand{\proof}{\noindent {\em Proof\ }: }
\def\endofproof{$\mbox{ }\hfill\Box$\par\vspace{1.8mm}\noindent}

\def\Rel{{\sf Rel}}
\def\o{^{\sf o}}

\def\Sh{{\sf Sh}}
\def\Sym{{\sf Sym}}
\def\Cat{{\sf Cat}}
\def\:{\colon}

\def\2{\textbf{2}}
\def\Set{{\sf Set}}
\def\op{^{\sf op}}
\def\dom{{\sf dom}}
\def\cod{{\sf cod}}
\def\Sup{{\sf Sup}}
\def\Dist{{\sf Dist}}
\def\Map{{\sf Map}}
\def\id{{\sf id}}
\def\cc{_{\sf cc}}
\def\sc{_{\sf sc}}
\def\s{_{\sf s}}
\def\C{{\cal C}}
\def\D{{\cal D}}
\def\E{{\cal E}}
\def\G{{\cal G}}
\def\K{{\cal K}}
\def\Q{{\cal Q}}
\def\R{{\cal R}}
\def\T{{\cal T}}
\def\bbA{\mathbb{A}}
\def\bbB{\mathbb{B}}
\def\bbC{\mathbb{C}}

\def\bbX{\mathbb{X}}
\def\tensor{\otimes}
\def\<{\langle}
\def\>{\rangle}
\def\eqref#1{(\ref{#1})}
\def\inv{^{-1}}

\title{Symmetry and Cauchy completion \\ of quantaloid-enriched categories}
\author{Hans Heymans\footnote{Department of Mathematics and Computer Science, University of Antwerp, Middelheimlaan 1, 2020 Antwerpen, Belgium, {\tt hans.heymans@ua.ac.be}}\ \ and Isar Stubbe\footnote{Laboratoire de Math\'ematiques Pures et Appliqu\'ees, Universit\'e du Littoral-C\^ote d'Opale, 50 rue F. Buisson, 62228 Calais, France, {\tt isar.stubbe@lmpa.univ-littoral.fr}}}
\date{June 23, 2011\footnote{Submitted May 4, 2010; accepted June 13, 2011; finalised June 23, 2011.}}

\begin{document}

\maketitle

\begin{abstract}
We formulate an elementary condition on an involutive quantaloid $\Q$ under which there is a distributive law from the Cauchy completion monad over the symmetrisation comonad on the category of $\Q$-enriched categories. For such quantaloids, which we call Cauchy-bilateral quantaloids, it follows that the Cauchy completion of any symmetric $\Q$-enriched category is again symmetric. Examples include Lawvere's quantale of non-negative real numbers and Walters' small quantaloids of closed cribles.
\end{abstract}

\section{Introduction}\label{A}

A quantaloid $\Q$ is a category enriched in the symmetric monoidal closed category $\Sup$ of complete lattices and supremum-preserving functions. Viewing $\Q$ as a bicategory, it is natural to study categories, functors and distributors enriched in $\Q$. If $\Q$ comes equipped with an involution, it makes sense to consider symmetric $\Q$-enriched categories. An important application of quantaloid-enriched categories was discovered by R.F.C.~Walters [1981, 1982]: he proved that the topos of sheaves on a small site $(\C,J)$ is equivalent to the category of {\em symmetric and Cauchy complete} categories enriched in a suitable ``small quantaloid of closed cribles'' $\R(\C,J)$. A decade earlier, F.W.~Lawvere [1973] had already pointed out that the category of generalised metric spaces and non-expansive maps is equivalent to the category of categories enriched in the quantale (that is, a one-object quantaloid) $([0,\infty],\bigwedge,+,0)$ of extended non-negative real numbers. This is a symmetric quantale, hence it is trivially involutive; and here too the {\em symmetric and Cauchy complete} $[0,\infty]$-enriched categories are important, if only to connect with the classical theory of metric spaces. Crucial in both examples is thus the use of categories enriched in an involutive quantaloid $\Q$ which are both symmetric and Cauchy complete. R. Betti and R.F.C.~Walters [1982] therefore raised the question ``whether the Cauchy completion of a symmetric [quantaloid-enriched] category is again symmetric''. That is to say, they ask whether it is possible to {\em restrict} the Cauchy completion functor $(-)\cc\:\Cat(\Q)\to\Cat(\Q)$ along the embedding $\Sym\Cat(\Q)\to\Cat(\Q)$ of symmetric $\Q$-categories. They show that the answer to their question is affirmative for both $\R(\C,J)$ and $[0,\infty]$, by giving an {\em ad hoc} proof in each case; they also give an example of an involutive quantale for which the answer to their question is negative. Thus, it depends on the base quantaloid $\Q$ whether or not the Cauchy completion of a symmetric $\Q$-category is again symmetric.

In this paper we address this issue in a slightly different manner to produce a single, simple argument for both Walters' small quantaloids of closed cribles and Lawvere's quantale of non-negative real numbers, thus giving perhaps a more decisive answer to Betti and Walters' question. The embedding $\Sym\Cat(\Q)\to\Cat(\Q)$ has a right adjoint $(-)\s\:\Cat(\Q)\to\Sym\Cat(\Q)$, which we call `symmetrisation'. We aim to {\em extend} the Cauchy completion functor along the symmetrisation functor. To wit, we define an obvious `symmetric completion' $(-)\sc\:\Sym\Cat(\Q)\to\Sym\Cat(\Q)$ (Proposition \ref{a3.0})
which, by construction, comes with a natural transformation
$$\xy
\xymatrix@=8ex{
\Cat(\Q)\ar[r]^{(-)\cc} & \Cat(\Q) \\
\Sym\Cat(\Q)\ar[r]_{(-)\sc}\ar[u]^{\mathrm{incl.}} & \Sym\Cat(\Q)\ar[u]_{\mathrm{incl.}}}
\POS(19,-7)\drop{K}\POS(17,-9)\drop{\rotatebox{145}{$\Longrightarrow$}}
\endxy$$
whose components $K_{\bbA}\:\bbA\sc\to\bbA\cc$ are full embeddings. Considering its mate
$$\xy
\xymatrix@=8ex{
\Cat(\Q)\ar[r]^{(-)\cc}\ar[d]_{(-)\s} & \Cat(\Q)\ar[d]^{(-)\s} \\
\Sym\Cat(\Q)\ar[r]_{(-)\sc} & \Sym\Cat(\Q)}
\POS(15,-7)\drop{L}\POS(17,-10)\drop{\rotatebox{35}{$\Longrightarrow$}}
\endxy$$
we formulate an elementary necessary-and-sufficient condition on $\Q$ under which $L$ is a natural isomorphism (Theorem \ref{a4}); in that case, we say that $\Q$ is a {\em Cauchy-bilateral} quantaloid (Definition \ref{XXX}). If $\Q$ is Cauchy-bilateral, then $K$ is in fact the identity transformation, thus in particular is the Cauchy completion of any symmetric $\Q$-category again symmetric. And moreover, as a corollary, we obtain a distributive law of the Cauchy completion monad over the symmetrisation comonad on $\Cat(\Q)$ (Corollary \ref{a7}). In a separate section we point out a number of examples, including Walters' small quantaloids of closed cribles and Lawvere's quantale of non-negative real numbers.

For an overview of the theory of quantaloid-enriched categories, and a list of appropriate historical references, we refer to [Stubbe, 2005], whose notations we adopt.

\section{Symmetric quantaloid-enriched categories}\label{S1}

In this section, after quickly recalling the notion of involutive quantaloid $\Q$, we give the obvious definition of symmetric $\Q$-category and explain how any $\Q$-category can be symmetrised. Examples are postponed to Section \ref{S3}.

\begin{definition}\label{a-1}
A {\em quantaloid} is a $\Sup$-enriched category. An {\em involution} on a quantaloid $\Q$ is a $\Sup$-functor $(-)\o\:\Q\op\to\Q$ which is the identity on objects and satisfies $f^{\sf oo}=f$ for any morphism $f$ in $\Q$. The pair $(\Q,(-)\o)$ is then said to form an {\em involutive quantaloid}.
\end{definition}
We shall often simply speak of ``an involutive quantaloid $\Q$'', leaving the notation for the involution understood. Note that the above definition is equivalent to an apparently weaker condition: in fact, any function $f\mapsto f\o$ on the morphisms of a quantaloid $\Q$ such that $f\leq g$ implies $f\o\leq g\o$, $(g\circ f)\o=f\o\circ g\o$, and $f^{\sf oo}=f$, is an involution. It is furthermore clear that an involution is an isomorphism between $\Q$ and $\Q\op$.

Whenever a morphism $f\:A\to B$ in a quantaloid (or in a locally ordered category, for that matter) is supposed to be a left adjoint, we write $f^*$ for its right adjoint. In many examples there is a big difference between the involute $f\o$ and the adjoint $f^*$ of a given morphism $f$, so morphisms for which involute and adjoint coincide, deserve a name:
\begin{definition}\label{a0}
In a quantaloid $\Q$ with involution $f\mapsto f\o$, an {\em ${\sf o}$-symmetric left adjoint} (or simply {\em symmetric left adjoint} if the context makes the involution clear) is a left adjoint whose right adjoint is its involute.
\end{definition}
Precisely as we write $\Map(\Q)$ for the category of left adjoints in $\Q$ (this notation being motivated by the widespread use of the word ``map'' synonymously with ``left adjoint''), we shall write $\Sym\Map(\Q)$ for the category of symmetric left adjoints. 

Recall that a category $\bbA$ enriched in a quantaloid $\Q$ consists of a set $\bbA_0$ of objects, each $x\in\bbA_0$ having a type $ta\in\Q_0$, and for any $x,y\in\bbA_0$ there is a hom-arrow $\bbA(y,x)\:tx\to ty$ in $\Q$, subject to associativity and unit requirements: $\bbA(z,y)\circ\bbA(y,x)\leq\bbA(z,x)$ and $1_{tx}\leq\bbA(x,x)$ for all $x,y,z\in\bbA_0$. A functor $F\:\bbA\to\bbB$ between such $\Q$-categories is an object-map $x\mapsto Fx$ such that $tx=t(Fx)$ and $\bbA(y,x)\leq\bbB(Fy,Fx)$ for all $x,y\in\bbA$. Such a functor is smaller than a functor $G\:\bbA\to\bbB$ if $1_{tx}\leq\bbB(Fx,Gx)$ for every $x\in\bbA$. With obvious composition one gets a locally ordered 2-category $\Cat(\Q)$ of $\Q$-categories and functors.

For two objects $x,y\in\bbA$, the hom-arrows $\bbA(y,x)$ and $\bbA(x,y)$ thus go in opposite directions. Hence, to formulate a notion of ``symmetry'' for $\Q$-categories, it is far too strong to require $\bbA(y,x)=\bbA(x,y)$. Instead, at least for involutive quantaloids, we better do as follows [Betti and Walters, 1982]:
\begin{definition}\label{a1}
Let $\Q$ be a small involutive quantaloid, with involution $f\mapsto f\o$. A $\Q$-category $\bbA$ is {\em ${\sf o}$-symmetric} (or {\em symmetric} if there is no confusion about the involved involution) when $\bbA(x,y)=\bbA(y,x)\o$ for every two objects $x,y\in\bbA$.
\end{definition}
We shall write $\Sym\Cat(\Q)$ for the full sub-2-category of $\Cat(\Q)$ determined by the symmetric $\Q$-categories; it is easy to see that the local order in $\Sym\Cat(\Q)$ is in fact symmetric (but not anti-symmetric). The full embedding $\Sym\Cat(\Q)\hookrightarrow\Cat(\Q)$ has a right adjoint functor\footnote{But the right adjoint is not a 2-functor, for this would imply the local order in $\Cat(\Q)$ to be symmetric, so this is not a 2-adjunction.}:
\begin{equation}\label{a2}
\Sym\Cat(\Q)\xymatrix@=8ex{\ar@{}[r]|{\perp}\ar@<1mm>@/^2mm/[r]^{\mathrm{incl.}} & \ar@<1mm>@/^2mm/[l]^{(-)\s}}\Cat(\Q).
\end{equation}
This `symmetrisation' sends a $\Q$-category $\bbA$ to the symmetric $\Q$-category $\bbA\s$ whose objects (and types) are those of $\bbA$, but for any two objects $x,y$ the hom-arrow is 
$$\bbA\s(y,x):=\bbA(y,x)\wedge\bbA(x,y)\o.$$ 
A functor $F\:\bbA\to\bbB$ is sent to $F\s\:\bbA\s\to\bbB\s\:a\mapsto Fa$. Quite obviously, the counit of this adjunction has components $S_{\bbA}\:\bbA\s\to\bbA\:a\mapsto a$.

Recall that a distributor $\Phi\:\bbA\dist\bbB$ between $\Q$-categories consists of arrows $\Phi(y,x)\:tx\to ty$ in $\Q$, one for each $(x,y)\in\bbA_0\times\bbB_0$, subject to action axioms: $\bbB(y',y)\circ\Phi(y,x)\leq\Phi(y',x)$ and $\Phi(y,x)\circ\bbA(x,x')\leq\Phi(y,x')$ for all $y,y'\in\bbB_0$ and $x,x'\in\bbA_0$. The composite of such a distributor with another $\Psi\:\bbB\dist\bbC$ is written as $\Psi\tensor\Phi\:\bbA\dist\bbC$, and its elements are computed with a ``matrix formula'': for $x\in\bbA_0$ and $z\in\bbC_0$,
$$(\Psi\tensor\Phi)(z,x)=\bigvee_{y\in\bbB_0}\Psi(z,y)\circ\Phi(y,x).$$
Parallel distributors can be compared elementwise, and in fact one gets a (large) quantaloid $
\Dist(\Q)$ of $\Q$-categories and distributors. Each functor $F\:\bbA\to\bbB$ determines an adjoint pair of distributors: $\bbB(-,F-)\:\bbA\dist\bbB$, with elements $\bbB(y,Fx)$ for $(x,y)\in\bbA_0\times\bbB_0$, is left adjoint to $\bbB(F-,-)\:\bbB\dist\bbA$ in the quantaloid $\Dist(\Q)$. These distributors are said to be `represented by $F$'. (More generally, a (necessarily left adjoint) distributor $\Phi\:\bbA\dist\bbB$ is `representable' if there exists a (necessarily essentially unique) functor $F\:\bbA\to\bbB$ such that $\Phi=\bbB(-,F-)$.) This amounts to a 2-functor 
\begin{equation}\label{a2.0}
\Cat(\Q)\to\Map(\Dist(\Q))\:\Big(F\:\bbA\to\bbB\Big)\mapsto\Big(\bbB(-,F-)\:\bbA\dist\bbB\Big).
\end{equation}

We shall write $\Sym\Dist(\Q)$ for the full subquantaloid of $\Dist(\Q)$ determined by the symmetric $\Q$-categories. It is easily verified that the involution $f\mapsto f\o$ on the base quantaloid $\Q$ extends to the quantaloid $\Sym\Dist(\Q)$: explicitly, if $\Phi\:\bbA\dist\bbB$ is a distributor between symmetric $\Q$-categories, then so is $\Phi\o\:\bbB\dist\bbA$, with elements $\Phi\o(a,b):=\Phi(b,a)\o$. And if $F\:\bbA\to\bbB$ is a functor between symmetric $\Q$-categories, then the left adjoint distributor represented by $F$ has the particular feature that it is a symmetric left adjoint in $\Sym\Dist(\Q)$ (in the sense of Definition \ref{a0}). That is to say, the functor in \eqref{a2.0} restricts to the symmetric situation as
\begin{equation}\label{a2.0.0}
\Sym\Cat(\Q)\to\Sym\Map(\Sym\Dist(\Q))\:\Big(F\:\bbA\to\bbB\Big)\mapsto\Big(\bbB(-,F-)\:\bbA\dist\bbB\Big),
\end{equation}
obviously giving a commutative diagram
\begin{equation}\label{a2.1}
\begin{array}{c}
\xymatrix@=8ex{
\Cat(\Q)\ar[r] & \Map(\Dist(\Q)) \\
\Sym\Cat(\Q)\ar[r]\ar[u]^{\mathrm{incl.}} & \Sym\Map(\Sym\Dist(\Q))\ar[u]_{\mathrm{incl.}}}
\end{array}
\end{equation}

\section{Cauchy completion and symmetric completion}\label{S2}

A $\Q$-category $\bbA$ is said to be `Cauchy complete' when each left adjoint distributor with codomain $\bbA$ is represented by a functor [Lawvere, 1973], that is, when for each $\Q$-category $\bbB$ the functor in \eqref{a2.0} determines an equivalence 
$$\Cat(\Q)(\bbB,\bbA)\simeq\Map(\Dist(\Q))(\bbB,\bbA).$$ 
It is equivalent to require this only for left adjoint presheaves\footnote{A `presheaf' on $\bbA$ is a distributor into $\bbA$ whose domain is a one-object category with an identity hom-arrow. Writing $*_X$ for the one-object $\Q$-category whose single object $*$ has type $X\in\Q_0$ and whose single hom-arrow is the identity $1_X$, a presheaf is then typically written as $\phi\:*_X\dist\bbA$. (These are really the {\em contravariant} presheaves on $\bbA$; the {\em covariant} presheaves are the distributors from $\bbA$ to $*_X$. In this paper, however, we shall only consider contravariant presheaves.)} on $\bbA$. The full inclusion of the Cauchy complete $\Q$-categories in $\Cat(\Q)$ admits a left adjoint:
\begin{equation}\label{z1}
\Cat\cc(\Q)\xymatrix@=8ex{\ar@{}[r]|{\perp}\ar@<-1mm>@/_2mm/[r]_{\mathrm{full\ incl.}} & \ar@<-1mm>@/_2mm/[l]_{(-)\cc}}\Cat(\Q).
\end{equation}
That is to say, each $\Q$-category $\bbA$ has a Cauchy completion $\bbA\cc$, which can be computed explicitly as follows: objects are the left adjoint presheaves on $\bbA$, the type of such a left adjoint $\phi\:*_X\dist\bbA$ is $X\in\Q$, and for another such $\psi\:*_Y\dist\bbA$ the hom-arrow $\bbA\cc(\psi,\phi)\:X\to Y$ in $\Q$ is the single element of the composite distributor $\psi^*\tensor\phi$ (where $\psi\dashv\psi^*$). The component at $\bbA\in\Cat(\Q)$ of the unit of this adjunction is a suitable corestriction of the Yoneda embedding: $Y_{\bbA}\:\bbA\to\bbA\cc\:x\mapsto\bbA(-,x)$. It is straightforward that $(-)\cc\:\Cat(\Q)\to\Cat(\Q)$ sends a functor $F\:\bbA\to\bbB$ to $F\cc\:\bbA\cc\to\bbB\cc\:\phi\mapsto\bbB(-,F-)\tensor\phi$. (For details, see e.g.\ [Stubbe, 2005, Section 7].)

The Cauchy completion can of course be applied to a symmetric $\Q$-category, but the resulting Cauchy complete category need not be symmetric anymore (see Example \ref{e5})! That is to say, the functor $(-)\cc\:\Cat(\Q)\to\Cat(\Q)$ does not restrict to $\Sym\Cat(\Q)$ in general. However, its very definition suggests the following modification: 
\begin{definition}\label{z2}
Let $\Q$ be a small involutive quantaloid. A symmetric $\Q$-category $\bbA$ is {\em symmetrically complete} if, for any symmetric $\Q$-category $\bbB$, the functor in \eqref{a2.0.0} determines an equivalence 
$$\Sym\Cat(\Q)(\bbB,\bbA)\simeq\Sym\Map(\Sym\Dist(\Q))(\bbB,\bbA).$$
\end{definition}
In analogy with the notion of Cauchy completeness of a $\Q$-category, it is straightforward to check the following equivalent expressions:
\begin{proposition}\label{a5}
Let $\Q$ be a small involutive quantaloid. For a symmetric $\Q$-category $\bbA$, the following conditions are equivalent:
\begin{enumerate}
\item $\bbA$ is symmetrically complete,
\item for any symmetric $\Q$-category $\bbB$, every symmetric left adjoint distributor $\Phi\:\bbB\dist\bbA$ is representable,
\item for any $X\in\Q$, every symmetric left adjoint presheaf $\phi\:*_X\dist\bbA$ is representable.
\end{enumerate}
\end{proposition}
And precisely as the Cauchy completion of a $\Q$-category can be computed explicitly with left adjoint presheaves, we can do as follows for the symmetric completion of a symmetric $\Q$-category:
\begin{proposition}\label{a3.0}
Let $\Q$ be a small involutive quantaloid. The full embedding of the symmetrically complete symmetric $\Q$-categories in $\Sym\Cat(\Q)$ admits a left adjoint:
\begin{equation}\label{z3}
\Sym\Cat\sc(\Q)\xymatrix@=8ex{\ar@{}[r]|{\perp}\ar@<-1mm>@/_2mm/[r]_{\mathrm{full\ incl.}} & \ar@<-1mm>@/_2mm/[l]_{(-)\sc}}\Sym\Cat(\Q).
\end{equation}
Explicitly, for a symmetric $\Q$-category $\bbA$, its symmetric completion $\bbA\sc$ is the full subcategory of $\bbA\cc$ determined by the {\em symmetric} left adjoint presheaves. The component at $\bbA\in\Sym\Cat(\Q)$ of the unit of this adjunction is a corestriction of the Yoneda embedding: $Y_{\bbA}\:\bbA\to\bbA\sc\:x\mapsto\bbA(-,x)$.
\end{proposition}
Note that $(-)\sc\:\Sym\Cat(\Q)\to\Sym\Cat(\Q)$ sends a functor $F\:\bbA\to\bbB$ between symmetric $\Q$-categories to $F\sc\:\bbA\sc\to\bbB\sc\:\phi\mapsto\bbB(-,F-)\tensor\phi$. Indeed, because $\bbB$ is symmetric, the distributor $\bbB(-,F-)$ is a symmetric left adjoint, hence its composition with $\phi\in\bbA\sc$ gives an object of $\bbB\sc$.

All this now raises a natural question: given any $\Q$-category $\bbA$, how does the symmetrisation of its Cauchy completion relate to the symmetric completion of its symmetrisation? It is clear that there is a natural transformation
\begin{equation}\label{a3.3}
\begin{array}{c}
\xy
\xymatrix@=8ex{
\Cat(\Q)\ar[r]^{(-)\cc} & \Cat(\Q) \\
\Sym\Cat(\Q)\ar[r]_{(-)\sc}\ar[u]^{\mathrm{incl.}} & \Sym\Cat(\Q)\ar[u]_{\mathrm{incl.}}}
\POS(20,-7)\drop{K}\POS(17,-9)\drop{\rotatebox{145}{$\Longrightarrow$}}
\endxy
\end{array}
\end{equation}
whose components are the full embeddings $K_{\bbA}\:\bbA\sc\to\bbA\cc\:\phi\mapsto\phi$ of which the construction, in Proposition \ref{a3.0}, of the symmetric completion speaks. From the calculus of mates [Kelly and Street, 1974] at least part of the following statement is then straightforward:
\begin{proposition}\label{a3.1}
Let $\Q$ be an involutive quantaloid. There is a natural transformation 
\begin{equation}\label{a3.2}
\begin{array}{c}
\xy
\xymatrix@=8ex{
\Cat(\Q)\ar[r]^{(-)\cc}\ar[d]_{(-)\s} & \Cat(\Q)\ar[d]^{(-)\s} \\
\Sym\Cat(\Q)\ar[r]_{(-)\sc} & \Sym\Cat(\Q)}
\POS(15,-7)\drop{L}\POS(17,-10)\drop{\rotatebox{35}{$\Longrightarrow$}}
\endxy
\end{array}
\end{equation}
whose component at $\bbA$ in $\Cat(\Q)$ is the full embedding\footnote{Recall that $S_{\bbA}\:\bbA\s\to\bbA\:a\mapsto a$ is the counit of the adjunction in diagram \eqref{a2}. Trivial as it may seem, it plays a crucial role throughout this section.}
$$L_\bbA\:(\bbA\s)\sc\to(\bbA\cc)\s\:\phi\mapsto\bbA(-,S_\bbA-)\tensor\phi.$$
Moreover, each $(\bbA\cc)\s$ is symmetrically complete.
\end{proposition}
\proof
It is straightforward to check that the explicit definition of $L$ is indeed obtained as the mate of the natural transformation in diagram \eqref{a3.3}. 

Given $\phi,\psi\in(\bbA\s)\sc$, consider the distributors
$$\xymatrix@=7ex{
{*_X}\ar[r]|{\distsign}^{\phi} & {\bbA\s}\ar@<-1ex>[d]|{\distsign}_{\bbA(-,S_{\bbA}-)}\ar@(ul,ur)|{\distsign}^{\bbA\s}\ar[r]|{\distsign}^{\psi\o} & {*_Y} \\
 & \bbA\ar@<-1ex>[u]|{\distsign}_{\bbA(S_{\bbA}-,-)} }$$
of which we know that $\bbA\s=\bbA(S_{\bbA}-,S_{\bbA}-)\wedge\bbA(S_{\bbA}-,S_{\bbA}-)\o$, and compute that
\begin{eqnarray*}
\psi\o\tensor\phi
 & = & \psi\o\tensor\bbA\s\tensor\phi \\
 & = & \psi\o\tensor(\bbA(S_{\bbA}-,S_{\bbA}-)\wedge\bbA(S_{\bbA}-,S_{\bbA}-)\o)\tensor\phi \\
 & = & (\psi\o\tensor\bbA(S_{\bbA}-,S_{\bbA}-)\tensor\phi)\wedge(\psi\o\tensor\bbA(S_{\bbA}-,S_{\bbA}-)\o\tensor\phi) \\
 & = & (\psi\o\tensor\bbA(S_{\bbA}-,S_{\bbA}-)\tensor\phi)\wedge(\phi\o\tensor\bbA(S_{\bbA}-,S_{\bbA}-)\tensor\psi)\o \\
 & = & (\psi\o\tensor\bbA(S_{\bbA}-,-)\tensor\bbA(-,S_{\bbA}-)\tensor\phi)\wedge(\phi\o\tensor\bbA(S_{\bbA}-,-)\tensor\bbA(-,S_{\bbA}-)\tensor\psi)\o \\
 & = & (L_{\bbA}(\psi)^*\tensor L_{\bbA}(\phi))\wedge (L_{\bbA}(\phi)^*\tensor L_{\bbA}(\psi))\o
\end{eqnarray*}
which asserts precisely the fully faithfulness of $L_{\bbA}$. (To pass from the second to the third line, we use that $\psi\o\tensor-\tensor\phi$ preserves infima, due to $\psi\dashv\psi\o$ and $\phi\dashv\phi\o$. From the third to the fourth line we use the involution on $\Sym\Dist(\Q)$ provided by the involution on $\Q$. And from line four to line five we use that $\bbA(S_{\bbA}-,S_{\bbA}-)=\bbA(S_{\bbA}-,-)\tensor\bbA(-,S_{\bbA}-)$.)

For any $a\in(\bbA\s)_0$ it is straightforward that $L_{\bbA}(\bbA\s(-,a))=\bbA(-,S_{\bbA}a)$. Putting $\psi=\bbA\s(-,a)$ in the previous calculation, we thus find for any $\phi\in(\bbA\s)\sc$ that
$$\bbA\s(a,-)\tensor\phi=\Big(\bbA(S_{\bbA}a,-)\tensor L_{\bbA}(\phi)\Big)\wedge\Big(L_{\bbA}(\phi)^*\tensor\bbA(-,S_{\bbA}a)\Big)\o.$$
Letting $a$ vary in $\bbA\s$, this shows that 
\begin{equation}\label{a3.4}
\phi=\Big(\bbA(S_{\bbA}-,-)\tensor L_{\bbA}(\phi)\Big)\wedge\Big(L_{\bbA}(\phi)^*\tensor\bbA(-,S_{\bbA}-)\Big)\o
\end{equation}
which implies that $L_{\bbA}$ is injective on objects.

As for the final part of the proposition, suppose that $\bbC$ is a Cauchy complete $\Q$-category and that $\phi\:*_X\dist\bbC\s$ is a symmetric left adjoint. Then there exists a $c\in\bbC_0$ such that $L_\bbC(\phi)=\bbC(-,c)$ and we use the formula in \eqref{a3.4} to compute that 
$$\phi=\Big(\bbC(S_{\bbC}-,-)\tensor L_\bbC(\phi)\Big)\wedge\Big(L_\bbC(\phi)^*\tensor\bbC(-,S_{\bbC}-)\Big)\o=\bbC(S_{\bbC}-,c)\wedge\bbC(c,S_{\bbC}-)\o=\bbC\s(-,c).$$ 
Therefore $\bbC\s$ is symmetrically complete. This of course applies to $\bbA\cc$.
\endofproof

Whereas the previous proposition establishes a {\it comparison} between $(\bbA\s)\sc$ and $(\bbA\cc)\s$, we shall now study when these two constructions {\it coincide}. This is related with the symmetrisation not only of $\Q$-categories and functors, but also of left adjoint distributors. We start by putting the formula in \eqref{a3.4} in a broader context:
\begin{lemma}\label{a3.5}
If $\Psi\:\bbA\dist\bbB$ is a left adjoint distributor between categories enriched in a small involutive quantaloid $\Q$, then the distributor
$$\Psi\s:=\Big(\bbB(S_{\bbB}-,-)\tensor\Psi\tensor\bbA(-,S_{\bbA}-)\Big)\wedge\Big(\bbA(S_{\bbA}-,-)\tensor\Psi^*\tensor\bbB(-,S_{\bbB}-)\Big)\o\:\bbA\s\dist\bbB\s$$
satisfies $\Psi\s\tensor(\Psi\s)\o\leq\bbB\s$. Therefore $\Psi\s$ is a symmetric left adjoint if and only if $\bbA\s\leq(\Psi\s)\o\tensor\Psi\s$; and if this is the case then it follows that $\Psi=\bbB(-,S_{\bbB}-)\tensor\Psi\s\tensor\bbA(S_{\bbA}-,-)$. 
\end{lemma}
\proof
It is clear that $\Psi\s\:\bbA\s\dist\bbB\s$ is a distributor: it is the infimum of two distributors, the first term of which is a composite of three distributors, and the second term is the involute of a composite of three distributors (which makes sense because domain and codomain of this composite are symmetric $\Q$-categories). Precisely because $\Psi\s$ is a distributor between symmetric $\Q$-categories, it makes sense to speak of its involute $(\Psi\s)\o$, and it is straightforward to compute that
\begin{eqnarray*}
\Psi\s\tensor(\Psi\s)\o
 & \leq & \Big(\bbB(S_{\bbB}-,-)\tensor\Psi\tensor\bbA(-,S_{\bbA}-)\Big)\tensor\Big(\bbA(S_{\bbA}-,-)\tensor\Psi^*\tensor\bbB(-,S_{\bbB}-)\Big)^{\sf oo} \\
 & \leq & \bbB(S_{\bbB}-,-)\tensor\Psi\tensor\bbA(-,-)\tensor\Psi^*\tensor\bbB(-,S_{\bbB}-) \\
 & \leq & \bbB(S_{\bbB}-,-)\tensor\bbB(-,-)\tensor\bbB(-,S_{\bbB}-) \\
 & = & \bbB(S_{\bbB}-,S_{\bbB}-)
\end{eqnarray*}
and therefore, by involution, also $\Psi\s\tensor(\Psi\s)\o\leq\bbB(S_{\bbB}-,S_{\bbB}-)\o$ holds, from which we can conclude that $\Psi\s\tensor(\Psi\s)\o\leq\bbB(S_{\bbB}-,S_{\bbB}-)\wedge\bbB(S_{\bbB}-,S_{\bbB}-)\o=\bbB\s(-,-)$ as claimed.

Now $\Psi\s$ is a symmetric left adjoint if and only if $\Psi\s\dashv(\Psi\s)\o$, and because the counit inequality of this adjunction always holds, this adjunction is equivalent to the truth of the unit inquality $\bbA\s\leq(\Psi\s)\o\tensor\Psi\s$. Suppose that this is indeed the case, then we can compute that
\begin{eqnarray*}
\Psi
 & = & \Psi\tensor\bbA(-,S_{\bbA}-)\tensor\bbA\s(-,-)\tensor\bbA(S_{\bbA}-,-) \\
 & \leq & \Psi\tensor\bbA(-,S_{\bbA}-)\tensor(\Psi\s)\o\tensor\Psi\s\tensor\bbA(S_{\bbA}-,-) \\
 & \leq & \Psi\tensor\bbA(-,S_{\bbA}-)\tensor\Big(\bbA(S_{\bbA}-,-)\tensor\Psi^*\tensor\bbB(-,S_{\bbB}-)\Big)^{\sf oo}\tensor\Psi\s\tensor\bbA(S_{\bbA}-,-) \\
 & \leq & \Psi\tensor\bbA(-,-)\tensor\Psi^*\tensor\bbB(-,S_{\bbB}-)\tensor\Psi\s\tensor\bbA(S_{\bbA}-,-) \\
 & \leq & \bbB\tensor\bbB(-,S_{\bbB}-)\tensor\Psi\s\tensor\bbA(S_{\bbA}-,-) \\
 & \leq & \bbB(-,S_{\bbB}-)\tensor\Big(\bbB(S_{\bbB}-,-)\tensor\Psi\tensor\bbA(-,S_{\bbA}-)\Big)\tensor\bbA(S_{\bbA}-,-) \\
 & \leq & \bbB(-,-)\tensor\Psi\tensor\bbA(-,-) \\
 & = & \Psi
\end{eqnarray*}
which means that $\Psi=\bbB(-,S_{\bbB}-)\tensor\Psi\s\tensor\bbA(S_{\bbA}-,-)$ as claimed.
\endofproof
The notation introduced in the previous lemma will be used in the remainder of this section. In particular shall we use it in the next proposition.
\begin{proposition}\label{a3.6}
For a category $\bbA$ enriched in a small involutive quantaloid $\Q$, the following conditions are equivalent:
\begin{enumerate}
\item\label{x1} the functor $L_{\bbA}\:(\bbA\s)\sc\to(\bbA\cc)\s$ from Proposition \ref{a3.1} is surjective on objects (and therefore an isomorphism, with inverse $\psi\mapsto\psi\s$),
\item\label{x2} for every left adjoint presheaf $\psi\:*_X\dist\bbA$, the presheaf $\psi\s\:*_X\dist\bbA\s$ is a symmetric left adjoint,
\item\label{x3} for every left adjoint distributor $\Psi\:\bbX\dist\bbA$, the distributor $\Psi\s\:\bbX\s\dist\bbA\s$ is a symmetric left adjoint.
\end{enumerate}
\end{proposition}
\proof
$(\ref{x1}\Leftrightarrow\ref{x2})$ From (the proof of) Proposition \ref{a3.1} we know that $L_{\bbA}$ is injective on objects, and that $\phi=(L_{\bbA}(\phi))\s$ for any symmetric left adjoint $\phi\:*_X\dist\bbA\s$ (this is the formula in \eqref{a3.4} rewritten with the notation introduced in Lemma \ref{a3.5}, taking into account that the domain of $\phi$ is the symmetric $\Q$-category $*_X$, so that $S_{*_X}$ is the identity functor on $*_X$). To say that $L_{\bbA}$ is surjective on objects thus means that for any left adjoint $\psi\:*_X\dist\bbA$ there exists a (necessarily unique) symmetric left adjoint $\phi\:*_X\dist\bbA\s$ such that $L_{\bbA}(\phi)=\psi$. Thus indeed $\psi\s=\phi$ is a symmetric left adjoint. Conversely, if we assume that for every left adjoint presheaf $\psi\:*_X\dist\bbA$ the presheaf $\psi\s\:*_X\dist\bbA\s$ is a symmetric left adjoint, then Lemma \ref{a3.5} implies $L(\psi\s)=\psi$ so that $L_{\bbA}$ is surjective on objects.

$(\ref{x3}\Leftrightarrow\ref{x2})$ One implication is trivial. For the other, by Lemma \ref{a3.5} we only need to prove that 
$$\bbX\s(y,x)\leq(\Psi\s)\o(y,-)\tensor\Psi\s(-,x)$$ 
for every left adjoint $\Psi\:\bbX\dist\bbA$ and every $x,y\in(\bbX\s)_0$. But for every $x\in\bbX_0$ we have a left adjoint presheaf $\Psi(-,x)\:*_{tx}\dist\bbA$ and by hypothesis thus also a symmetric left adjoint presheaf $\Psi(-,x)\s\:*_{tx}\dist\bbA\s$. Because $\Psi(-,x)\s=\Psi\s(-,x)$ and $(\Psi\s)\o(y,-)=(\Psi\s(-,y))\o=(\Psi(-,y)\s)\o$, the sought-after inequation is equivalent to
$$\bbX\s(y,x)\leq(\Psi(-,y)\s)\o\tensor\Psi(-,x)\s.$$
Using the adjunction $\Psi(-,x)\s\dashv(\Psi(-,x)\s)\o$ this is in turn equivalent to 
$$ \bbX\s(y,x)\tensor(\Psi(-,x)\s)\o\leq(\Psi(-,y)\s)\o$$
which is an instance of the action inequality $\bbX\s\tensor(\Psi\s)\o\leq(\Psi\s)\o$ for  $(\Psi\s)\o\:\bbA\s\dist\bbX\s$.
\endofproof
Now we have everything in place to prove our main theorem, establishing in particular an elementary necessary-and-sufficient condition on the base quantaloid $\Q$ under which $(\bbA\s)\sc\cong(\bbA\cc)\s$ holds for {\em every} $\Q$-category $\bbA$.
\begin{theorem}\label{a4}
For a small involutive quantaloid $\Q$, the following conditions are equivalent:
\begin{enumerate}
\item\label{t2} each functor $L_{\bbA}\:(\bbA\s)\sc\to(\bbA\cc)\s$ as in Proposition \ref{a3.1} is an isomorphism (making diagram \eqref{a3.2} commute up to isomorphism), 
\item\label{t1} for every left adjoint presheaf $\psi\:*_X\dist\bbA$, the presheaf $\psi\s\:*_X\dist\bbA\s$ is a symmetric left adjoint,
\item\label{t3} for each left adjoint distributor $\Psi\:\bbA\dist\bbB$, the distributor $\Psi\s\:\bbA\s\dist\bbB\s$ is a symmetric left adjoint,
\item\label{z5} the inclusion $\Sym\Map(\Sym\Dist(\Q))\to\Map(\Dist(\Q))$ admits a right adjoint making the following square commute:
$$\xymatrix@=8ex{
\Cat(\Q)\ar[r]\ar[d]_{(-)\s} & \Map(\Dist(\Q))\ar@{.>}[d] \\
\Sym\Cat(\Q)\ar[r] & \Sym\Map(\Sym\Dist(\Q))}$$
\item\label{t4} for each family $(f_i\:X\to X_i, g_i\:X_i\to X)_{i\in I}$ of morphisms in $\Q$,
$$\left.\begin{array}{c}
\forall j,k\in I:\ f_k\circ g_j\circ f_j\leq f_k \\[1ex]
\forall j,k\in I:\ g_j\circ f_j\circ g_k\leq g_k \\[1ex]
1_X\leq\displaystyle\bigvee_{i\in I}g_i\circ f_i
\end{array}\right\}\Longrightarrow\
1_X\leq\bigvee_{i\in I}(g_i\wedge f_i\o)\circ(g_i\o\wedge f_i).$$
\end{enumerate}
In fact, the right adjoint of which the fourth statement speaks, is
\begin{equation}\label{z6}
(-)\s\:\Map(\Dist(\Q))\to\Sym\Map(\Sym\Dist(\Q))\:\Big(\Psi\:\bbA\dist\bbB\Big)\mapsto\Big(\Psi\s\:\bbA\s\dist\bbB\s\Big).
\end{equation}
\end{theorem}
\proof 
$(\ref{t2}\Leftrightarrow\ref{t1}\Leftrightarrow\ref{t3})$ Are taken care of in Proposition \ref{a3.6}.

$(\ref{t3}\Rightarrow\ref{z5})$ With the help of Lemma \ref{a3.5}, it can be checked that the left adjoint distributor $\bbC(-,S_{\bbC}-)\:\bbC\s\dist\bbC$ displays $\bbC\s$ as the co\-re\-flec\-tion of a $\Q$-category $\bbC$ along the inclusion $\Sym\Map(\Sym\Dist(\Q))\to\Map(\Dist(\Q))$: if $\bbA$ is a symmetric $\Q$-category and $\Psi\:\bbA\dist\bbC$ is a left adjoint distributor, then by assumption we have that $\Psi\s\:\bbA\dist\bbC\s$ is a symmetric left adjoint distributor such that $\Psi=\bbC(-,S_{\bbC}-)\tensor\Psi\s$; and if $\Phi\:\bbA\dist\bbC\s$ would be another symmetric left adjoint distributor such that $\Psi=\bbC(-,S_{\bbC}-)\tensor\Phi$ (and therefore also $\Psi^*=\Phi\o\tensor\bbC(S_{\bbC}-,-)$), then necessarily 
\begin{eqnarray*}
\Psi\s
 & = & \Big(\bbC(S_{\bbC}-,-)\tensor\Psi\Big)\wedge\Big(\Psi^*\tensor\bbC(-,S_{\bbC}-)\Big)\o \\
 & = & \Big(\bbC(S_{\bbC}-,-)\tensor\bbC(-,S_{\bbC}-)\tensor\Phi\Big)\wedge\Big(\Phi\o\tensor\bbC(S_{\bbC}-,-)\tensor\bbC(-,S_{\bbC}-)\Big)\o \\
 & = & \Big(\bbC(S_{\bbC}-,S_{\bbC}-)\tensor\Phi\Big)\wedge\Big(\bbC(S_{\bbC}-,S_{\bbC}-)\o\tensor\Phi\Big) \\
 & = & \Big(\bbC(S_{\bbC}-,S_{\bbC}-)\wedge\bbC(S_{\bbC}-,S_{\bbC}-)\o\Big)\tensor\Phi \\
 & = & \bbC\s\tensor\Phi \\
 & = & \Phi
\end{eqnarray*}
(To pass from the third to the fourth line we use that $-\tensor\Phi$ preserves infima because $\Phi$ is a left adjoint.) By general theory for adjoint functors, these coreflections $\bbC(-,S_{\bbC}-)\:\bbC\s\dist\bbC$ determine the right adjoint, which is the functor given in \eqref{z6}.

$(\ref{z5}\Rightarrow\ref{t2})$ Suppose that $G\:\Sym\Map(\Dist(\Q))\to\Map(\Dist(\Q))$ is right adjoint to the inclusion and makes the square commute; thus $G$ necessarily acts on objects as $\bbC\mapsto\bbC\s$. Writing $\varepsilon_{\bbC}\:\bbC\s\dist\bbC$ for the counit of the adjunction, $G(\varepsilon_{\bbC})$ is necessarily the identity distributor on $\bbC\s$. But on the other hand, by commutativity of the diagram, $G$ also sends $\bbC(-,S_{\bbC}-)\:\bbC\s\dist\bbC$ to the identity distributor on $\bbC\s$. By general theory for adjoint functors, $G(\varepsilon_{\bbC})=G(\bbC(-,S_{\bbC}-))$ implies $\varepsilon_{\bbC}=\bbC(-,S_{\bbC}-)$. Thus $\bbC(-,S_{\bbC}-)\:\bbC\s\dist\bbC$ enjoys a universal property, saying in particular that: for every left adjoint $\phi\:*_X\dist\bbC$ there is a unique symmetric left adjoint $\psi\:*_X\dist\bbC\s$ such that $\bbC(-,S_{\bbC}-)\tensor\psi=\phi$. In other words, $L_{\bbC}\:(\bbC\s)\sc\to(\bbC\cc)\s\:\psi\mapsto\bbC(-,S_{\bbC}-)\tensor\psi$ is surjective on objects, and therefore an isomorphism.

$(\ref{t1}\Rightarrow\ref{t4})$ Putting $\bbA_0=I$, $ti=X_i$ and $\bbA(j,i)=f_j\circ g_j\vee\delta_{ij}$ defines a $\Q$-category $\bbA$ (the ``Kronecker delta'' $\delta_{ij}\:X_i\to X_j$ denotes the identity when $i=j$ and the zero morphism otherwise), and putting $\psi(i)=f_i$ defines a presheaf $\psi\:*_X\dist\bbA$ with a right adjoint $\psi^*\:\bbA\dist *_X$ which is given by $\psi^*(i)=g_i$. By hypothesis we infer that $\psi\s\:*_X\dist\bbA\s$ is a symmetric left adjoint. This means in particular that $1_X\leq(\psi\s)\o\tensor\psi\s$, or in other terms $1_X\leq\bigvee_i(f_i\o\wedge g_i)\circ(f_i\wedge g_i\o)$, as wanted.

$(\ref{t4}\Rightarrow\ref{t1})$ If $\psi\:*_X\dist\bbA$ is a left adjoint, the family $(\psi(a)\:X\to ta,\psi^*(a)\:ta\to X)_{a\in\bbA_0}$ of morphisms in $\Q$ is easily seen to satisfy the conditions in the hypothesis, thus
$$1_X\leq\bigvee_{a\in\bbA_0}\Big(\psi(a)\o\wedge\psi^*(a)\Big)\circ\Big(\psi(a)\wedge\psi^*(a)\Big)\o.$$ 
The right hand side is exactly $\psi\s\o\tensor\psi\s$ so this is equivalent to $\psi\s\dashv(\psi\s)\o$ (cf.\ Lemma \ref{a3.5}). 
\endofproof
For further reference, we give a name to those quantaloids (small or large) that satisfy the fifth condition in the above Theorem \ref{a4}:
\begin{definition}\label{XXX}
A quantaloid $\Q$ is {\em Cauchy-bilateral} if it is involutive (with involution $f\mapsto f\o$) and  for each family $(f_i\:X\to X_i, g_i\:X_i\to X)_{i\in I}$ of morphisms in $\Q$,
$$\left.\begin{array}{c}
\forall j,k\in I:\ f_k\circ g_j\circ f_j\leq f_k \\[1ex]
\forall j,k\in I:\ g_j\circ f_j\circ g_k\leq g_k \\[1ex]
1_X\leq\displaystyle\bigvee_{i\in I}g_i\circ f_i
\end{array}\right\}\Longrightarrow\
1_X\leq\bigvee_{i\in I}(g_i\wedge f_i\o)\circ(g_i\o\wedge f_i).$$
\end{definition}
Thus, a {\em small} Cauchy-bilateral quantaloid $\Q$ is precisely one that satisfies the equivalent conditions in Theorem \ref{a4}.

To finish this section we explain an important consequence of Theorem \ref{a4}, containing an answer to R. Betti and R.F.C.~Walters' [1982] question about the symmetry of the Cauchy completion of a symmetric category:
\begin{corollary}\label{a7}
If $\Q$ is a small Cauchy-bilateral quantaloid, then
the following diagrams commute:
$$\xymatrix@=8ex{
\Cat(\Q)\ar[r]^{(-)\cc} & \Cat(\Q) \\
\Sym\Cat(\Q)\ar[r]_{(-)\cc}\ar[u]^{\mathrm{incl.}} & \Sym\Cat(\Q)\ar[u]_{\mathrm{incl.}}}
\hspace{5ex}
\xymatrix@=8ex{
\Cat(\Q)\ar[r]^{(-)\s} & \Cat(\Q) \\
\Cat\cc(\Q)\ar[r]_{(-)\s}\ar[u]^{\mathrm{incl.}} & \Cat\cc(\Q)\ar[u]_{\mathrm{incl.}}}
$$
\end{corollary}
\proof
Suppose that the equivalent conditions in Theorem \ref{a4} hold. If $\psi\:*_X\dist\bbA$ is a left adjoint presheaf on a symmetric $\Q$-category, then $\psi\s\:*_X\dist\bbA$ is a symmetric left adjoint presheaf which satisfies $\psi=\bbA(-,S_{\bbA}-)\tensor\psi\s=\psi\s$ (by Lemma \ref{a3.5} and symmetry of $\bbA$). So $\psi$ is necessarily a symmetric left adjoint. Hence the full embedding $\bbA\sc\hookrightarrow\bbA\cc$ is surjective-on-objects, or in other words, $\bbA\sc=\bbA\cc$. Therefore the Cauchy completion of a symmetric $\Q$-category is symmetric, making the first square commute. Now suppose that $\bbC$ is a Cauchy complete category. In (the proof of) Proposition \ref{a3.1} it was stipulated that $\bbC\s$ is symmetrically complete, so -- knowing now that the symmetric completion and the Cauchy completion of any symmetric $\Q$-category coincide -- it follows that $\bbC\s$ is Cauchy complete too, making the second square commute.
\endofproof
This corollary implies that, whenever a small $\Q$ is Cauchy-bilateral, there is a {\em distributive law} [Beck, 1969; Street, 1972; Power and Watanabe, 2002] of the Cauchy completion monad over the symmetrisation comonad on the category $\Cat(\Q)$. More precisely, the reflective subcategory $\Cat\cc(\Q)\to\Cat(\Q)$ is the category of algebras of a monad $(\T,\mu,\eta)$ on $\Cat(\Q)$; and similarly, the coreflective subcategory $\Sym\Cat(\Q)\to\Cat(\Q)$ is the category of coalgebras of a comonad $(\D,\delta,\varepsilon)$ on $\Cat(\Q)$. (The functors $\T$ and $\D$ are precisely $(-)\cc$ and $(-)\s$, of course.) As shown in [Power and Watanabe, 2002, Theorems 3.10 and 5.10], the commutativity of the squares in Corollary \ref{a7} is {\em equivalent} to the existence of a natural transformation $\lambda\:\T\circ\D\tto\D\circ\T$ making the following diagrams commute:
$$\xymatrix@=8ex{
\T\T\D\bbC\ar[r]^{\mu_{\D\bbC}}\ar[d]_{\T\lambda_{\bbC}} & \T\D\bbC\ar[dd]^{\lambda_{\bbC}} \\
\T\D\T\bbC\ar[d]_{\lambda_{\T\bbC}} & & \D\bbC\ar[ul]_{\eta_{\D\bbC}}\ar[dl]^{\D\eta_{\bbC}} \\
\D\T\T\bbC\ar[r]_{\D\mu_{\bbC}} & \D\T\bbC}
\hspace{3ex}
\xymatrix@=8ex{
 & \T\D\bbC\ar[dl]_{\T\varepsilon_{\bbC}}\ar[r]^{\T\delta_{\bbC}}\ar[dd]_{\lambda_{\bbC}} & \T\D\D\bbC\ar[d]^{\lambda_{\D\bbC}} \\
\T\bbC & & \D\T\D\bbC\ar[d]^{\D\lambda_{\bbC}} \\
 & \D\T\bbC\ar[ul]^{\varepsilon_{\T\bbC}}\ar[r]_{\delta_{\T\bbC}} & \D\D\T\bbC}$$
This, in turn, says exactly that $\lambda$ is a {\em distributive law of the monad $\T$ over the comonad $\D$} [Power and Watanabe, 2002, Definition 6.1]. Because $\T$ and $\D$ arise from (co)reflective subcategories, there is {\em at most one} such distributive law; its components are necessarily 
$$\lambda_{\bbC}\:(\bbC\s)\cc\to(\bbC\cc)\s\:\phi\mapsto\bbC(-,S_{\bbC}-)\tensor\phi.$$ 
(Thus $\lambda_{\bbC}$ is precisely the functor $L_{\bbC}$ of Proposition \ref{a3.1}, reckoning that -- under the conditions of Theorem \ref{a4} -- the symmetric completion of a symmetric $\Q$-category coincides with its Cauchy completion.) It is a consequence of the general theory of distributive laws that the monad $\T$ restricts to the category of $\D$-coalgebras, that the comonad $\D$ restricts to the category of $\T$-algebras, and that the categories of (co)algebras for these restricted (co)monads are equivalent to each other and are further equivalent to the category of so-called {\em $\lambda$-bialgebras} [Power and Watanabe, 2002, Corollary 6.8]. In the case at hand, a $\lambda$-bialgebra is simply a $\Q$-category which is both symmetric and Cauchy-complete (the ``$\lambda$-compatibility'' between algebra and coalgebra structure is trivially satisfied), and a morphism between $\lambda$-bialgebras is simply a functor between such $\Q$-categories.

\section{Examples}\label{S3}

\begin{example}[Commutative quantales]\label{e1}
A quantale is, by definition, a one-object quantaloid. (For some authors, a quantale need not be unital, so for them it is {\em not} a one-object quantaloid; but for us, a quantale is always unital.) Put differently, a quantale is a monoid in the monoidal category $\Sup$ (whereas a quantaloid is a $\Sup$-enriched category). Obviously, a quantale $Q$ is commutative if and only if the identity function $1_Q\:Q\to Q$ is an involution.
\end{example}

As we shall point out below, many an interesting involutive quantaloid $\Q$ satisfies the following condition: 
\begin{definition}\label{YYY}
A quantaloid $\Q$ is {\em strongly Cauchy-bilateral} when it is involutive (with involution $f\mapsto f\o$) and for any family $(f_i\:X\to X_i, g_i\:X_i\to X)_{i\in I}$ of morphisms in $\Q$,
$$1_X\leq\bigvee_ig_i\circ f_i \ \Longrightarrow \ 1_X\leq\bigvee_i(f_i\o\wedge g_i)\circ(f_i\wedge g_i\o).$$
\end{definition}
Obviously, a strongly Cauchy-bilateral $\Q$ is Cauchy-bilateral in the sense of Definition \ref{XXX}. For a so-called integral quantaloid -- that is, when the top element of each $\Q(X,X)$ is $1_X$ -- Cauchy-bilaterality and strong Cauchy-bilaterality are equivalent notions, but in general the latter is {\em strictly stronger} than the former:
\begin{example}\label{e7}
We take an example of an involutive quantale from [Resende, 2007, Example 3.18]: the complete lattice with Hesse diagram
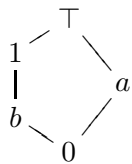
$$\xymatrix@R=0ex@C=1.5ex{
 & \top & \\
1\ar@{-}[ur] \\
 & & a\ar@{-}[uul] \\
b\ar@{-}[uu] \\
 & 0\ar@{-}[ul]\ar@{-}[uur]}$$
together with the commutative multiplication defined by $a\circ\top=\top$, $a\circ a=b$ and $a\circ b=a$. This gives a quantale which we can equip with the identity involution. It is straightforward to check that this quantale is Cauchy-bilateral, but not strongly so.
\end{example}

\begin{example}[Generalised metric spaces]\label{e3}
The condition for strong Cauchy-bilaterality is satisfied by the integral and commutative quantale $Q=\left([0,\infty],\bigwedge,+,0\right)$ with its trivial involution: for any family $(a_i,b_i)_{i\in I}$ of pairs of elements of $[0,\infty]$, if $\bigwedge_i(a_i+b_i)\leq 0$ is assumed then 
$$\bigwedge_i\left(\max\{a_i,b_i\}+\max\{a_i,b_i\}\right)=2\cdot\bigwedge_i\max\{a_i,b_i\}\leq 2\cdot\bigwedge_i(a_i+b_i)\leq 0.$$
This ``explains'' the well known fact that the Cauchy completion of a symmetric generalised metric space [Lawvere, 1973] is again symmetric.
\end{example}

\begin{example}[Locales]\label{e2}
Any locale $(L,\bigvee,\wedge,\top)$ is a commutative (hence trivially involutive) and integral quantale; it is easily checked that $L$ is strongly Cauchy-bilateral. Splitting the idempotents of the $\Sup$-monoid $(L,\wedge,\top)$ gives an integral quantaloid with an obvious involution; it too is strongly Cauchy-bilateral.
\end{example}

\begin{example}[Groupoid-quantaloids with canonical involution]\label{e4}
For a category $\C$, let $\Q(\C)$ be the quantaloid with the same objects as $\C$ but where $\Q(\C)(X,Y)$ is the complete lattice of subsets of $\C(X,Y)$, composition is done ``pointwise'' (for $S\subseteq\C(X,Y)$ and $T\subseteq\C(Y,Z)$ let $T\circ S:=\{t\circ s\mid t\in T,s\in S\}$) and the identity on an object $X$ is the singleton $\{1_X\}$. With local suprema in $\Q(\C)$ given by union, it is straightforward to check that $\Q(\C)$ is a quantaloid; it is the free quantaloid on the category $\C$.

If $\G$ is a groupoid, then $\Q(\G)$ comes with a {\em canonical involution} $S\mapsto S\o:=\{s\inv\mid s\in S\}$. For any family $(T_i\subseteq\G(X,X_i),S_i\subseteq\G(X_i,X))_{i\in I}$ we can prove that
$$1_X\in\bigcup_iS_i\circ T_i\ \Longrightarrow\ 1_X\in\bigcup_i(S_i\o\cap T_i) \circ(S_i\cap T_i\o).$$
Indeed, the premise says that there is a $i_0\in I$ for which we have $x\in S_{i_0}$ and $y\in T_{i_0}$ such that $1_X=x \circ y$ in $\G$. But then $y\in S_{i_0}\o\cap T_{i_0}$ and $x\in S_{i_0}\cap T_{i_0}\o$, so that the conclusion follows. That is to say, $\Q(\G)$ is strongly Cauchy-bilateral.
\end{example}

\begin{example}[Commutative group-quantales with trivial involution]\label{e5}
For a commutative group $(G,\cdot,1)$, also the group-quantale $\Q(G)$ is commutative, and -- in contrast with the above example -- it can therefore be equipped with the {\em trivial involution} $S\mapsto S\o:=S$. Betti and Walters [1982] gave a simple example of such a commutative group-quantale with trivial involution for which the Cauchy completion of a symmetric enriched category is not necessarily symmetric. We repeat it here: Let $G=\{1,a,b\}$ be the commutative group defined by $a\cdot a=b$, $b\cdot b=a$ and $a\cdot b=1$; then Betti and Walters showed that the Cauchy completion of the (symmetric) singleton $\Q(G)$-category whose hom is $\{1\}$, is not symmetric. In fact, the pair $(\{a\},\{b\})$ of elements of $\Q(G)$ does satisfy the premise but not the conclusion of the condition in Definition \ref{XXX}: thus, in retrospect, Theorem \ref{a4} predicts that there must exist a symmetric category whose Cauchy completion is no longer symmetric.
\end{example}

There is a common generalisation of Examples \ref{e2} and \ref{e4}, due to [Walters, 1982]: given a small site $(\C,J)$, there is an involutive quantaloid $\R(\C,J)$ such that the category of symmetric and Cauchy complete $\R(\C,J)$-categories is equivalent to $\Sh(\C,J)$. We shall spell out this important example, and show that it is strongly Cauchy-bilateral (and thus also satisfies the equivalent conditions in Theorem \ref{a4}). In retrospect, this proves that the {\em symmetric and Cauchy complete $\R(\C,J)$-categories} can be computed as the {\em Cauchy completions of the symmetric $\R(\C,J)$-categories}.

\begin{example}[Quantaloids determined by small sites]\label{e2+4}
If $\C$ is a small category, then the small quantaloid $\R(\C)$ of {\em cribles} in $\C$ is the full sub-quantaloid of $\Rel(\Set^{\C\op})$ whose objects are the representable presheaves. It is useful to have an explicit description. We write a {\em span} in $\C$ as $(f,g)\:D\span C$, and intend it to be a pair of arrows with $\dom(f)=\dom(g)$, $\cod(f)=C$ and $\cod(g)=D$. A crible $R\:D\crible C$ is then a set of spans $D\span C$ such that for any $(f,g)\in R$ and any $h\in\C$ with $\cod(h)=\dom(f)$, also $(f\circ h,g\circ h)\in R$. Composition in $\R(\C)$ is obvious: for $R\:D\crible C$ and $S\:E\crible D$ the elements of $R\circ S\:E\crible C$ are the spans $(f,g):E\span C$ for which there exists a morphism $h\in\C$ such that $(f,h)\in R$ and $(h,g)\in S$. The identity crible $\id_C\:C\crible C$ is the set $\{(f,f)\mid\cod(f)=C\}$, and the supremum of a set of cribles is their union. In fact, $\R(\C)$ is an involutive quantaloid: the involute $R\o\:C\crible D$ of a crible $R\:D\crible C$ is obtained by reversing the spans in $R$. 

If $J$ is a Grothendieck topology on the category $\C$, then there is a locally left exact nucleus\footnote{A nucleus $j$ on a quantaloid $\Q$ is a lax functor $j\:\Q\to\Q$ which is the identity on objects and such that each $j\:\Q(X,Y)\to\Q(X,Y)$ is a closure operator; it is locally left exact if it preserves finite infima of arrows. If $j\:\Q\to\Q$ is a nucleus on a quantaloid, then there is a quotient quantaloid $\Q_j$ of $j$-closed morphisms, that is, those $f\in\Q$ for which $j(f)=f$.} $j\:\R(\C)\to\R(\C)$ sending a crible $R\:D\crible C$ to
$$j(R):=\{(f,g)\:D\span C\mid\exists S\in J(\dom(f)):\forall s\in S,\ (g\circ s,f\circ s)\in R\}.$$
Conversely, if $j\:\R(\C)\to\R(\C)$ is a locally left exact nucleus, then
$$J(C):=\{S\mbox{ is a sieve on }C\mid\id_C\leq j(\{(s,s)\mid ,s\in S\})\}$$
defines a Grothendieck topology on $\C$. These procedures are each other's inverse [Betti and Carboni, 1983; Rosenthal, 1996]. For a small site $(\C,J)$ we write, following [Walters, 1982], $\R(\C,J)$ for the quotient quantaloid $\R(\C)_j$ where $j$ is the nucleus determined by the Grothendieck topology $J$; it is the {\em small quantaloid of closed cribles} determined by the site $(\C,J)$. Because every locally left exact nucleus on $\R(\C)$ preserves the involution, $\R(\C,J)$ too is involutive. (Walters [1982] originally called $\R(\C,J)$ a `bicategory of relations', wrote it as $\Rel(\C,J)$, and called its arrows `relations'. To avoid confusion with other constructions that have been called `bicategories of relations' since then, we prefer to speak of `small quantaloids of closed cribles'. For an axiomatic study of these, we refer to [Heymans and Stubbe, 2011].)

Any locale $L$ can be thought of as a site $(\C,J)$, where $\C$ is the ordered set $L$ and $J$ is its so-called canonical topology (so $J(u)$ is the set of all covering families of $u\in L$): it is easily verified that $\R(\C,J)$ is then isomorphic (as involutive quantaloid) to the quantaloid obtained by splitting the idempotents in the $\Sup$-monoid $L$. And if $\G$ is a small groupoid and $J$ is the smallest Grothendieck topology on $\G$, then the quantaloid of relations $\R(\G,J)$ equals the quantaloid of cribles $\R(\G)$, which in turn is isomorphic (as involutive quantaloid) to the free quantaloid $\Q(\G)$ with its canonical involution. Indeed, any crible $R\:X\crible Y$ in $\G$ determines the subset $F(R):=\{h\inv\circ g\mid (g,h)\in R\}$ of $\G(X,Y)$. Conversely, for any subset $S$ of $\G(X,Y)$ let $G(S)$ be the smallest crible containing the set of spans $\{(1_X,s)\mid s\in S\}$ in $\G$. Then $R\mapsto F(R)$ and $S\mapsto G(S)$ extend to functors $F\:\R(\G)\to\Q(\G)$ and $G\:\Q(\G)\to\R(\G)$ which are each other's inverse and which preserve the involution. Hence both Examples \ref{e2} and \ref{e4} are covered by the construction of the quantaloid $\R(\C,J)$ from a small site $(\C,J)$. 

Now we show that $\R(\C,J)$ is strongly Cauchy-bilateral, as in Definition \ref{YYY}. This claim is equivalent to saying that for any family $(F_i\:X\crible X_i,G_i\:X_i\crible X)_{i\in I}$ of cribles in $\C$ we have
$$j(\id_X)\subseteq j\Big(\bigcup_i j\Big(j(G_i)\circ j(F_i)\Big)\Big)\ \Longrightarrow\ j(\id_X)\subseteq j\Big(\bigcup_i j\Big((j(F_i)\o\cap j(G_i)\circ (j(F_i)\cap j(G_i)\o)\Big)\Big)$$
in the involutive quantaloid $\R(\C)$ with left exact nucleus $j$ constructed from $J$.

Take the left-hand side: it is equivalent (by general computations with the nucleus $j$) to $\id_X\subseteq j(\bigcup_i G_i\circ F_i)$
in $\R(\C)$. By definition this means that, for any morphism $x$ in $\C$ with $\cod(x)=X$, $(x,x)\in j(\bigcup_i G_i\circ F_i)$, but because $j(\bigcup_i G_i\circ F_i)$ is a crible, this is equivalent to requiring simply that $(1_X,1_X)\in j(\bigcup_i G_i\circ F_i)$ (and here $1_X$ denotes the identity on $X$ in $\C$). Spelling out the definition of the nucleus $j$ in terms of the Grothendieck topology $J$ this means that: there exists a covering sieve $S\in J(X)$ such that for all $s\in S$ there exists an $i_s\in I$ satisfying $(s,s)\in G_{i_s}\circ F_{i_s}$. In a similar fashion the right-hand side can be seen to say precisely that: there exists a covering sieve $S\in J(X)$ such that for all $s\in S$ there exists an $i_s\in I$ satisfying $(s,s)\in (F_{i_s}\o\cap G_{i_s})\circ(F_{i_s}\cap G_{i_s}\o)$. The implication is now straightforward.
\end{example}

The following example further generalises the previous one.

\begin{example}[Locally localic and modular quantaloids]\label{a30}
Following [Freyd and Scedrov, 1990] we say that a quantaloid $\Q$ is locally localic when each $\Q(X,Y)$ is a locale; and $\Q$ is modular if it is involutive and when for any morphisms $f\:Z\to Y, g\:Y\to X$ and $h\:Z\to X$ in $\Q$ we have $gf\wedge h\leq g(f\wedge g\o h)$ (or equivalently, $gf\wedge h\leq (g\wedge hf\o)f$). (Here we write the composition in $\Q$ by juxtaposition to avoid overly bracketed expressions.) In fact, every locally localic and modular quantaloid $\Q$ is strongly Cauchy-bilateral: suppose that  $(f_i\:X\to X_i, g_i\:X_i\to X)_{i\in I}$ is a family of morphisms in $\Q$ such that $1_X\leq\bigvee_ig_if_i$, then we can compute that:
\begin{eqnarray*}
1_X
 & = & 1_X\wedge\bigvee_i g_i f_i \\
 & = & \bigvee_i (1_X\wedge g_i f_i) \\
 & = & \bigvee_i (1_X\wedge (1_X\wedge g_i f_i)) \\
 & \leq & \bigvee_i(1_X\wedge g_i(g_i\o 1_X\wedge f_i)) \\
 & = & \bigvee_i(1_X\wedge g_i(g_i\o\wedge f_i)) \\
 & \leq & \bigvee_i(1_X(g_i\o\wedge f_i)\o\wedge g_i)(g_i\o\wedge f_i) \\
 & = & \bigvee_i (g_i\wedge f_i\o)(g_i\o\wedge f_i).
\end{eqnarray*}
To pass from the first to the second line we used that $\Q(X,X)$ is a locale, and both inequalities were introduced by use of the modular law.

Any small quantaloid of relations $\R(\C,J)$ is in fact locally localic and modular, and thus its strong Cauchy-bilaterality follows from the above computation. But to prove that $\R(\C,J)$ is modular, is not simpler than to prove directly that it satisfies the condition in Definition \ref{YYY}, as we did in Example \ref{e2+4}. (The quantale in Example \ref{e3} is locally localic but not modular, and the quantale in Example \ref{e7} is neither locally localic nor modular; but both are Cauchy-bilateral.)
\end{example}

\begin{example}[Sets and relations]\label{e8}
The quantaloid $\Rel$ of sets and relations is not small, but it is involutive (the involute of a relation is its opposite: $R\o=\{(y,x)\mid(x,y)\in R\}$) and it is strongly Cauchy-bilateral. In fact, this holds for any quantaloid $\Rel(\E)$ of internal relations in a Grothendieck topos $\E$, because it is modular and locally localic [Freyd and Scedrov, 1990].
\end{example}
There is a subtle difference between Examples \ref{e2+4} and \ref{e8}: the former deals with the {\em small} quantaloid $\R(\C,J)$ built from a small site, the latter deals with the {\em large} quantaloid $\Rel(\Sh(\C,J))$ of relations between the sheaves on that site. However, both constructions give rise to a Cauchy-bilateral quantaloid. We shall further analyse the interplay between these quantaloids in a forthcoming paper.

Finally we mention a difference between ``symmetric'' and ``discrete'' $\Q$-categories.
\begin{example}[Symmetric vs.\ discrete]\label{e9}
In any locally ordered category $\K$, an object $D$ is said to be {\em discrete} when, for any other object $X\in\K$, the order $\K(X,D)$ is symmetric. It is straightforward to verify that, whenever $\Q$ is a small Cauchy-bilateral quantaloid, every symmetric and Cauchy complete $\Q$-category is a discrete object of $\Cat\cc(\Q)$. However, not all discrete objects of $\Cat\cc(\Q)$ need to be symmetric, not even when $\Q$ is Cauchy-bilateral! A counterexample can be found in the theory of generalised metric spaces: Suppose that $X$ is a set and $R\subseteq X\times X$. For $x,y\in X$, a {\em path} from $x$ to $y$ is a sequence $\alpha=(x_0,x_1,\ldots,x_n)$ of elements of $X$ with $x_0=x$, $x_n=y$ and $(x_i,x_{i+1})\in R$ for all $i<n$; the {\em length} $l(\alpha)$ of such a path $\alpha$ is then $n$. For every $x\in X$, $(x)$ is a path from $x$ to $x$ of length $0$. It is easy to verify that 
$$d_R(x,y):=\bigwedge\{l(\alpha)\mid\alpha\mbox{ is a path from $x$ to $y$}\}$$ 
turns $X$ into a generalised metric space. (This infimum is a minimum, except when there is no path from $x$ to $y$, in which case $d_R(x,y)=\infty$.) Any such space $(X,d_R)$ is Cauchy complete, as is every generalised metric space with values in $\mathbb{N}\cup\{\infty\}$. And, in fact, it is discrete in the sense given above, because $x\leq y$ if and only if $0\geq d_R(x,y)$, so there is a path with length 0 from $x$ to $y$, which means that $x=y$. However, choosing $X=\{0,1\}$ and $R=\{(0,1)\}$ gives a non-symmetric example: $d_R(0,1)=1\neq\infty=d_R(1,0)$.
\end{example}

\end{document}